\documentclass[12pt]{article}
\usepackage{amsfonts,amsmath}
\def \beq {\begin{eqnarray}}
\def \eeq {\end{eqnarray}}
\def \beqn {\begin{eqnarray*}}
\def \eeqn {\end{eqnarray*}}

\newcommand{\halmos}{\rule{1ex}{1.4ex}}

\newcounter{for}[section]

\newtheorem{itlemma}{Lemma}[section]
\newtheorem{itproposition}[itlemma]{Proposition}
\newtheorem{theorem}[itlemma]{Theorem}
\newtheorem{itcorollary}[itlemma]{Corollary}
\newtheorem{itremark}[itlemma]{Remark}
\newtheorem{itremarks}[itlemma]{Remarks}
\newtheorem{itdefinition}[itlemma]{Definition}
\newtheorem{itexample}[itlemma]{Example}

\newenvironment{fact}{\begin{itfact}\rm}{\end{itfact}}
\newenvironment{claim}{\begin{itclaim}\rm}{\end{itclaim}}
\newenvironment{lemma}{\begin{itlemma}}{\end{itlemma}}
\newenvironment{remark}{\begin{itremark}\rm}{\end{itremark}}
\newenvironment{remarks}{\begin{itremarks} \rm}{\end{itremarks}}
\newenvironment{corollary}{\begin{itcorollary}}{\end{itcorollary}}
\newenvironment{proposition}{\begin{itproposition}}{\end{itproposition}}
\newenvironment{definition}{\begin{itdefinition}\rm}{\end{itdefinition}}
\newenvironment{example}{\begin{itexample}\rm}{\end{itexample}}
\newenvironment{proof}{\noindent {\em Proof}.\ \
}{\hspace*{\fill}$\halmos$\medskip}
\newcommand{\be}[1]{\addtocounter{for}{1} \begin{equation}\label{#1}}
\newcommand{\ee}{\end{equation}}
\newcommand{\bl}[1]{\begin{lemma}\label{#1}}
\newcommand{\br}[1]{\begin{remark}\label{#1}}
\newcommand{\brs}[1]{\begin{remarks}\label{#1}}
\newcommand{\bt}[1]{\begin{theorem}\label{#1}}
\newcommand{\bd}[1]{\begin{definition}\label{#1}}
\newcommand{\bp}[1]{\begin{proposition}\label{#1}}
\newcommand{\bc}[1]{\begin{corollary}\label{#1}}
\newcommand{\bfact}[1]{\begin{fact}\label{#1}}
\newcommand{\bex}[1]{\begin{example}\label{#1}}
\newcommand{\ec}{\end{corollary}}
\newcommand{\efact}{\end{fact}}
\newcommand{\eex}{\end{example}}
\newcommand{\el}{\end{lemma}}
\newcommand{\er}{\end{remark}}
\newcommand{\ers}{\end{remarks}}
\newcommand{\et}{\end{theorem}}
\newcommand{\ed}{\end{definition}}
\newcommand{\ep}{\end{proposition}}
\newcommand{\epr}{\end{proof}}
\newcommand{\bpr}{\begin{proof}}
\newcommand{\bcl}[1]{\begin{claim}\label{#1}}
\newcommand{\ecl}{\end{claim}}
\newcommand{\ecs}{\end{corollary}}
\newcommand{\eers}{\end{exercise}}
\newcommand{\eexs}{\end{example}}
\newcommand{\eems}{\end{example}}
\newcommand{\els}{\end{lemma}}
\newcommand{\eles}{\end{lemmaex}}
\newcommand{\ets}{\end{theorem}}
\newcommand{\eds}{\end{definition}}
\newcommand{\eps}{\end{proposition}}

\newcommand{\bi}{\begin{itemize}}
\newcommand{\ei}{\end{itemize}}
\newcommand{\ben}{\begin{enumerate}}
\newcommand{\een}{\end{enumerate}}

\def\vbar{\mathchoice{\vrule height6.3ptdepth-.5ptwidth.8pt\kern-.8pt}
   {\vrule height6.3ptdepth-.5ptwidth.8pt\kern-.8pt}
   {\vrule height4.1ptdepth-.35ptwidth.6pt\kern-.6pt}
   {\vrule height3.1ptdepth-.25ptwidth.5pt\kern-.5pt}}
\def\fudge{\mathchoice{}{}{\mkern.5mu}{\mkern.8mu}}
\def\bbc#1#2{{\rm \mkern#2mu\vbar\mkern-#2mu#1}}
\def\bbb#1{{\rm I\mkern-3.5mu #1}}
\def\bba#1#2{{\rm #1\mkern-#2mu\fudge #1}}
\def\bb#1{{\count4=`#1 \advance\count4by-64 \ifcase\count4\or\bba A{11.5}\or
   \bbb B\or\bbc C{5}\or\bbb D\or\bbb E\or\bbb F \or\bbc G{5}\or\bbb H\or
   \bbb I\or\bbc J{3}\or\bbb K\or\bbb L \or\bbb M\or\bbb N\or\bbc O{5} \or
   \bbb P\or\bbc Q{5}\or\bbb R\or\bbc S{4.2}\or\bba T{10.5}\or\bbc U{5}\or
   \bba V{12}\or\bba W{16.5}\or\bba X{11}\or\bba Y{11.7}\or\bba Z{7.5}\fi}}

\def \uchi{\underline{\chi}}

\def \P{{\cal{P}}}
\def \LL {{\cal{L}}}

\def \A {{\cal{A}}}

\def \G {{\cal{G}}}

\def \C {{\cal{C}}}
\def \D {{\cal{D}}}
\def \Db{\underline{{\mathcal D}}_b}
\def \Dh{{\bar{\mathcal D}}_b}

\newcommand{\ba}[1]{\addtocounter{for}{1} \begin{eqnarray}\label{#1}}
\newcommand{\ea}{\end{eqnarray}}

\def\sqr#1#2{{\vcenter{\vbox{\hrule height .#2pt
                             \hbox{\vrule width .#2pt height#1pt \kern#1pt
                                   \vrule width .#2pt}
                             \hrule height .#2pt}}}}

\def\pmb#1{\setbox0=\hbox{#1}%
   \kern-.025em\copy0\kern-\wd0
   \kern.05em\copy0\kern-\wd0
 \kern-.025em\raise.0433em\box0 }
\def\sqr#1#2{{\vcenter{\vbox{\hrule height.#2pt
     \hbox{\vrule width.#2pt height#1pt \kern#1pt
   \vrule width.#2pt}\hrule height.#2pt}}}}

\def \RR {{\mathcal R}}
\def\N{{\mathbb N}}
\def\Z{{\mathbb Z}}

\def\PP{{\mathbb P}}

\def\reff#1{(\ref{#1})}

\newcommand {\cro}[1] {\left[ {#1} \right]}
\newcommand {\acc}[1] {\left\{ {#1} \right\}}
\newcommand {\pare}[1] {\left( {#1} \right)}

\begin{document}

\title{Large deviations estimates for self-intersection local times for simple
random walk in $\Z^3$.}
\author{Amine Asselah \\C.M.I., Universit\'e de Provence,\\
39 Rue Joliot-Curie, \\F-13453 Marseille cedex 13, France\\
asselah@cmi.univ-mrs.fr}
\date{}
\maketitle
\begin{abstract}
We obtain large deviations estimates for the self-intersection
local times for a symmetric random walk in dimension 3.
Also, we show that the main contribution to making the self-intersection
large, in a time period of length $n$, comes from sites visited less than some
power of $\log(n)$.
This is opposite to the situation in dimensions larger or equal to 5.
Finally, we present two applications of our estimates: (i) to 
moderate deviations estimates for the range of a random walk, and (ii)
to moderate deviations for random walk in random sceneries. 
\end{abstract}

{\em Keywords and phrases}: self-intersection local times, random walk,
random sceneries.

{\em AMS 2000 subject classification numbers}: 60K35, 82C22,
60J25.

{\em Running head}: Large deviations for self-intersections in $d=3$.

\section{Introduction} 
In this paper, we focus on large deviations estimates for the self-intersection
local times (SILT) for a simple random walk in dimension 3. Thus,
$P_x$ denotes the law of a nearest neighbors
symmetric random walk $\{S_k,k\geq 0\}$ on $\Z^d$ starting at site
$x\in \Z^d$, and for any $y\in \Z^d$ and $n\in \N$, the local time $l_n(y)$ 
is the number of visits of $y$ up to time $n$.
The SILT process is denoted $\{\Sigma_n^2, n\in \N\}$ with
\be{def-SILT}
\Sigma_n^2=\sum_{x\in \Z^d}l^2_n(x).
\ee

This paper is a sequel to recent works \cite{ACa} and \cite{ACb} dealing with
dimensions $d\ge 5$. There, the initial motivation came
from establishing large deviations estimates for random walk in random sceneries:
in \cite{ACa}, this problem was reduced to estimating the distribution
of the size of the level sets of the local times. In other words, for
large $L$ and $t$, one needed to estimate the probability
of $\A_n(t,L):=\{|\{x: l_n(x)\sim t\}|>L\}$, where for a set $\Lambda$, we
denote by $|\Lambda|$ its cardinal. A key tool in \cite{ACa} was the following 
simple observation (see Lemma 2.1 of \cite{ACa}): when $d\ge 3$,
there is a constant $\kappa_d$ such that for any subset $\Lambda$
in $\Z^d$, 
\be{main-first}
P_0\pare{l_n(\Lambda)>t}\le \exp\pare{-\kappa_d \frac{t}{|\Lambda|^{2/d}}},
\quad\text{where}\quad l_n(\Lambda)=\sum_{x\in \Lambda} l_n(x).
\ee
Then, in order to use \reff{main-first}, $\A_n(t,L)$ was partitioned as follows
\[
\A_n(t,L)\subset \bigcup \acc{ l_n(\Lambda)>t |\Lambda|
:\ \Lambda\subset [-n,n]^d \quad\&\quad |\Lambda|=L}.
\]
Thus, the uniform estimate \reff{main-first} yielded
\be{ACa-1}
P_0\pare{\A_n(t,L)}\le 
\C_n(L) \exp\pare{-\kappa_d t L^{1-2/d}},\quad\text{with}\quad
\C_n(L):=|\acc{\Lambda\subset [-n,n]^d: |\Lambda|=L}| 
\ee
In \cite{ACa}, the combinatorial term $\C_n(L)$ in \reff{ACa-1}
had an innocuous r\^ole since $P_0\pare{\A_n(t_n,L)}$ was needed for a sequence
$\{t_n\}$ so large that the trivial bound $L<n/t_n$ made $\C_n(L)$ negligible
compared to $\exp(\kappa_d t_n L^{1-2/d})$. 
However, in \cite{ACb}, the combinatorial term ruined the naive bound \reff{ACa-1}.
Thus, $\A_n(t,L)$ was first transformed into a question for the SILT
\be{ACb-1}
\A_n(t,L)\subset \acc{ \sum_{x\in \Z^d} 1_{\{l_n(x)\sim t\}}l_n^2(x)\ge L t^2}.
\ee
Then the key estimate of \cite{ACb} (Lemma 2.1) relied on bounding the
self-intersection times {\it of a given level set of the local times} by the
intersection times for two independent half-trajectories over a {\it larger
level set}.
This observation which twisted an idea of Le Gall~\cite{legall} reads {\it morally} as
\be{main-second}
\acc{ \sum_{x:l_n(x)\sim t} l_n^2(x)\ge L t^2}
\quad\longleftrightarrow\quad
\acc{ \sum_{x:l_{n/2}(x)\le t} l_{n/2}(x)\tilde l_{n/2}(x)\ge L t^2}
\subset \acc{\tilde l_{n/2}\pare{\D}\ge L t},
\ee
where $\D:=\acc{l_{n/2}(x)\le t}$, and $\{\tilde l_n(x),x\in \Z^d\}$ is an
independent copy of the local times with law $\tilde P_0$. Thus, one reformulates
the key tool~\reff{main-first} in order 
to get rid of the combinatorial factor as follows:
\be{lem-ACa}
P_0\otimes\tilde P_0\pare{ \tilde l_{n/2}(\D)>z,\ |\D|<y}\le
\exp\pare{-\kappa_d F(z,y)},\quad\text{with}\quad
F(z,y)=\frac{z}{y^{2/d}}.
\ee
Thus, we can evaluate $P_0(\tilde l_{n/2}\pare{\D}\ge L t)$ with 
the help of \reff{lem-ACa} as soon as a good bound on $|\D|$ obtains.

Note however, that in \reff{main-second}, the sum on the right hand side
is over $\acc{x:l_{n/2}(x)\le t}$. This poses no trouble in $d\ge 5$, since
the main contribution comes from {\it large} level sets. However, this approach fails
in $d=3$ and $d=4$, and can at best bring a spurious logarithmic term as in
the upper bound of \reff{UB-prop1}. Besides, no indication can be extracted
as to which level set gives a dominant contribution.

In this paper, the approach is somewhat opposite to that of \cite{ACb}: we deal directly
with level sets' distribution which in turn provides new estimates
for the SILT process. The key idea is to transform any given level set of the
local times of $\{S_0,S_1,\dots,S_{2n}\}$ into two sets:
\begin{itemize}
\item The sites that at least one of the {\it half}
trajectories $\{S_n-S_{n-1},\dots,S_n-S_0\}$,
or $\{S_n-S_{n+1},\dots,S_n-S_{2n}\}$ visits {\it nearly as often} as the whole
trajectory.
\item The sites that both trajectories visit {\it enough} times.
\end{itemize}
Then, we iterate this procedure chopping each trajectories near its midpoint from
which stems two independent trajectories, and so 
forth until no piece of trajectories remains.
This seemingly innocent strategy
allows us to obtain some informations in dimension 3.

We show that the main contribution in making $\Sigma_n^2$ large
comes from sites which are ``not too often'' visited. This is drastically different
from the situation in $d\ge 5$, where only a few sites, where
$l_n(x)\sim \sqrt{\Sigma_n^2}$, contributed to
making $\Sigma_n^2$ large (see \cite{ACb}). In dimension 4, it is still an open problem
to understand which level sets give a dominant contribution to realize the
large deviation $\{\Sigma_n^2> ny\}$.
\bp{prop1} In dimension $d=3$, there are positive constants $\underline c,
\bar c$ such that for $y$ large enough
\be{UB-prop2bis}
\exp\pare{ -\underline c y^{2/3}n^{1/3} }\le
P_0\pare{ \sum_{x\in \Z^d} l_n^2(x)> ny}\le \exp\pare{ -\bar c y^{1/3} n^{1/3}}.
\ee
Moreover, there is $\uchi>0$ such that if $\bar\D:=
\acc{x:\ l_n(x)> \log(n)^{\uchi}}$, then there is $\epsilon>0$ such that
\be{UB-prop2}
P_0\pare{ \sum_{x\in \bar\D} l_n^2(x)> ny}\le \exp\pare{ -\bar c n^{1/3}
\log(n)^{\epsilon} }.
\ee
\ep
\br{rem-intro} It is a simple application of
Lemma 2.1 of \cite{ACb}, and of our moment computations in Lemma~\ref{lem-app2}
to obtain that in dimensions 3 and 4, for $y$ large enough
there are positive constants $\underline c, \bar c, \chi$ such that 
\be{UB-prop1}
\exp\pare{ -\underline c n^{1-2/d} }\le
P_0\pare{ \sum_{x\in \Z^d} l_n^2(x)> ny}\le \exp\pare{ -\bar c \frac{n^{1-2/d}}
{\log(n)^{\chi}} }.
\ee
However, the upper bound of \reff{UB-prop2bis} and most importantly
\reff{UB-prop2} require a new treatment of the level sets.
\er
A heuristic understanding of Proposition~\ref{prop1} comes from the following scenario
realizing the lower bound in~\reff{UB-prop1}: we localize the walk
in a ball $B(r_n)$ of radius $r_n$ with $r_n^3\sim n/y$. Indeed, assume that sites
of $B(r_n)$ are visited uniformly: for $x\in B(r_n)$, $l_n(x)\sim n/r_n^3
\sim y$, and thus $\Sigma_n^2\sim ny$. Now, the probability of staying in $B(r_n)$
a period of time $n$ is larger than $\exp(-C n/r_n^2)$ (for some $C>0$), 
which yields the right exponent.
However, we cannot say if, in the optimal strategy, the walk spends a fraction of its time
outside $B(r_n)$, as expected by the result 
of van den Berg, Bolthausen \& den Hollander\cite{bolt} concerning
the volume of the Wiener sausage, which is is the continuous counterpart 
of the range of the walk $\RR_n:=\acc{x: l_n(x)>0}$.
Indeed, a connection between the two problems (already noticed in \cite{ACb})
is as follows: 
\be{range-silt}
\frac{n}{|\RR_n|}\le \frac{\sum_{x\in \Z^d} l_n^2(x)}{n}
\Longrightarrow 
P_0\pare{|\RR_n|<\frac{n}{y}}\le P_0\pare{\sum_{x\in \Z^d} l_n^2(x)>y n}.
\ee
Note that in $d\ge 5$, the results of \cite{ACb} show that the range does not shrink
when realizing $\{\Sigma_n^2>ny\}$, whereas in $d=3$, the cost of the two deviations
(i.e. small $|\RR_n|$ and large $\Sigma_n^2$) correspond to the same speed $n^{1/3}$, and
it would be interesting to know whether $\RR_n$ shrinks to produce $\{\Sigma_n^2>ny\}$.

In dimension 2, large and moderate deviation principles are established for the
SILT for Brownian motion in Bass \& Chen \cite{bc04}, and for stable processes in
Bass, Chen \& Rosen\cite{bcr05}. 
Also, moderate deviations for the SILT and for the range of planar random walks were
recently obtained by Bass, Chen \& Rosen respectively in~\cite{bcr05b} and~\cite{bcr06}. 
The approach of \cite{bc04,bcr05,bcr05b,bcr06}
lies ultimately on the Donsker-Varadhan large deviation principle for the
Brownian occupation measure~\cite{DV}, and might not be adequate when the 
dominant strategy to perform the large deviations is not a localization. Finally, for
the $d=1$ case, we refer the reader to  Mansmann \cite{Ma91}, 
and Chen \& Li \cite{chen-li}.

We now present two applications of our estimates on self-intersections.
First, knowing that a random walk stays a time $n$ in a ball $B(r_n)$ with
$n/r_n^3\gg 1$, we show that typically
a proportion of the sites of $B(r_n)$ are visited about $n/|B(r_n)|$.
Let $\sigma(r)$ be the first time the random walk exits the ball $B(r)$ of radius $r$.
Also, we use the common notation $a_n=O(b_n)$ meaning that for some constant $A>0$,
$|a_n|\le A |b_n|$.

\bp{prop-range} Let $\{r_n\}$ be a sequence going to infinity with $r_n^3=O(n)$.
When $\epsilon_0$ and $\delta_0$ are small enough, we have
\be{visits-1}
\lim_{n\to\infty} P_0\pare{ \left|\acc{x:l_n(x)>\delta_0 \frac{n}{|B(r_n)|}}\right|
\ge \epsilon_0 |B(r_n)|\quad\Big|\!\!\Big|\quad \sigma(r_n)>n}=1.
\ee
\ep
\br{rem-range} Proposition~\ref{prop-range} is based on the following estimate.
For $y$ large enough, the inequality \reff{range-silt} and the
upper bound in \reff{UB-prop2bis} imply that there is a constant $\kappa$
such that
\be{boltup}
P_0\pare{|\RR_n|<\frac{n}{y}}\le \exp\pare{-\kappa y^{1/3} n^{1/3}}.
\ee
This is weaker than the asymptotics of van de Berg, Bolthausen \&
den Hollander\cite{bolt} for the volume of the Wiener sausage, and the proof is simpler. 
Also, to establish a lower bound similar to \reff{boltup}, 
note that the range is small if we localize the walk
in a ball $B(r_n)$ with $|B(r_n)|=n/y$. Thus,
\be{boltlow}
\{\sigma(r_n)>n\}\subset \{|\RR_n|<\frac{n}{y}\}\Longrightarrow
P_0\pare{|\RR_n|<\frac{n}{y}}\ge \exp\pare{-C\frac{n}{(n/y)^{2/3}}}=
e^{-Cy^{2/3} n^{1/3}}.
\ee
\er

Secondly, we establish moderate deviations estimates for random walk in random 
sceneries (RWRS), following the approach of \cite{ACb}. Thus, we consider a field
$\{\eta(x),x\in \Z^d\}$ independent of the random walk $\{S_k,n\in \N\}$, and made up of
centered i.i.d.\  with law denoted by $P_{\eta}$ and tail decay 
\be{eq-tail.1}
\lim_{t\to\infty} \frac{\log P_{\eta}(\eta(0)>t)}{t^{\alpha}}=-c,
\quad\text{for a positive constant  }c.
\ee
The random walk in random scenery is the process $X_n=\eta(S_0)+\dots+\eta(S_n)$.
We present asymptotics for the probability, averaged over both randomness,
that $\{X_n>n^{\beta}\}$ for $\beta>1/2$ and $\alpha\ge 1$ in dimension 3.
Our estimates are of the following type. For $\beta>1/2$, and $y$ large enough,
there are two positive constants $c_1,c_2$ such that if $\PP:=P_0\otimes
P_{\eta}$
\be{rwrs-type}
\exp\pare{-c_1 n^{\zeta}}\le \PP(X_n>y n^\beta)\le \exp\pare{-c_2 n^{\zeta}}.
\ee
Thus, the next result consists in characterizing the
exponent $\zeta$ as a function of $(\alpha,\beta)$.
\bp{prop-rwrs} Assume that dimension is 3.
\begin{itemize}
\item In region I$:=\{(\alpha,\beta):1\le\alpha,\  1/2<\beta\le 2/3\}$, 
we have $\zeta_{I}=2\beta-1$.
\item In region II$:=\{(\alpha,\beta):1\le\alpha<3/2,\ \beta>\frac{1+\alpha}{4-\alpha}\}$,
we have  $\zeta_{I\!I}=\beta\frac{\alpha}{1+\alpha}$.
\item In region III$:=\{(\alpha,\beta): 1\le\alpha,\ 
2/3<\beta<\min(1,\frac{1+\alpha}{(4-\alpha)^+})\}$, 
we have $\zeta_{I\!I\!I}=\frac{4}{5}\beta-\frac{1}{5}$.
\end{itemize}
\ep
\br{rem-rwrs} Compared with the situation in dimensions $d\ge 5$, we see
that region III, which corresponds to localizing the walk, has expanded
in $d=3$. Note also that the lower bounds in regions I and II are already
written in \cite{ACb}. Also, we refer to \cite{ACb} for a discussion of the
behaviour of the walk and the environment leading to the exponent $\zeta$ in
each region.
\er

Note that region IV$:=\{(\alpha,\beta): \alpha>2/3,\ \beta\ge 1\}$ is
treated in \cite{GKS}, where a large deviation principle is established.
Also, a regime with $\alpha<1$ is thoroughly studied in \cite{HGK}.

We prove Proposition~\ref{prop1} in Section~\ref{sec-main}, whose
Subsection~\ref{sec-new} is our main technical part. In Section~\ref{sec-low},
we establish Proposition~\ref{prop-range}, and the lower bound in Region III 
for Proposition~\ref{prop-rwrs}.
Finally, we have gathered in the Appendix a useful large deviation estimate
and moments computations for intersection local times in $d=3$ and $d=4$.

\section{Proof of Proposition~\ref{prop1}}\label{sec-main}
Note first that in order to obtain \reff{UB-prop2},
we do not need to worry about the contribution of
$\acc{x:l_n(x)\ge n^{1/3+\epsilon}}$, for $\epsilon>0$, since
in dimension $d=3$, $l_n(x)$ is bounded by a geometric variable and the
upper bound of \reff{UB-prop2} follows easily for 
$P(\acc{x:l_n(x)>n^{1/3+\epsilon}}\not=
\emptyset)$. Also, we set for simplicity $n=2^N$, and consider a subdivision
$\{N^{\alpha_j},j=1,\dots,M_N\}$ of $[1,2^{N(1/3+\epsilon)}]$,
with $\alpha_j=(j-1)\alpha+\uchi$,
for positive constants $\alpha,\uchi$ to be chosen later.
Note also that $M_N$ is of order $N/\log(N)$.
We now form the level sets of the local times
\be{eq-not.18}
\LL_j=\{x: N^{\alpha_j}\le l_{2^N}(x)< N^{\alpha_{j+1}}\}\quad\text{for }j>0,\quad
\text{and}\quad \LL_0=\{x: 1\le l_{2^N}(x)< N^{\uchi}\}.
\ee
Also, let $y_j=y/(2M_N)$ for $j>0$, and $y_0=y/2$ so that $y_0+\dots+y_{M_N}=y$.
We have the following decomposition 
\be{eq-not.19}
\begin{split}
\acc{ \sum_{x\in \Z^d}  l_{2^N}^2(x)>y2^N}& \subset\bigcup_{j=0}^{M_N}
\acc{ \sum_{\LL_j} l_{2^N}^2(x)>y_j2^N}\cup\acc{x:l_{2^N}(x)\ge 2^{N(\frac{1}{3}+
\epsilon)}}\\
\!\!\!&\subset\bigcup_{j=1}^{M_N}
\acc{ |\LL_j|>\frac{2^N y_j}{N^{2\alpha_{j+1}}} }
\cup \acc{ \sum_{\LL_0} l_{2^N}^2(x)>y_02^N }
\cup\acc{l_{2^N}\ge 2^{N(\frac{1}{3}+\epsilon)}}.
\end{split}
\ee
In Section~\ref{sec-new}, 
we deal with estimating the distribution of $|\LL_j|$ for $j>0$.
In Section~\ref{tilt}, we consider $\acc{\sum_{\LL_0} l_{2^N}^2(x)>n y_0}$.
Finally, in Section~\ref{lower} we repeat an argument of \cite{ACb}
to obtain the lower bound of Proposition~\ref{prop1}.
\subsection{Proof of \reff{UB-prop2}}\label{sec-new}
We relabel our original trajectory as $\{S_k^{(0)},k\in \N\}$ and
its local time as $\{l^{(0)}_{k,1},k\in \N\}$. We fix a time $2^N$, and build
from $\{S_0^{(0)},\dots,S_{2^N}^{(0)}\}$ 
two independent trajectories running for times $k\in \{0,\dots,2^{N-1}\}$
\be{eq-not.1}
S^{(1)}_{k,1}=S^{(0)}_{2^{N-1}}-S^{(0)}_{2^{N-1}-k} ,\quad\text{and}\quad
S^{(1)}_{k,2}=S^{(0)}_{2^{N-1}}-S^{(0)}_{2^{N-1}+k}.
\ee
We denote by $\{l^{(1)}_{2^{N-1},i}(x),x\in \Z^d\}$ the local times of
$\{S^{(1)}_{k,i}\}$ at time $2^{N-1}$ for $i=1,2$. Likewise, 
we proceed inductively,
and consider at {\it generation} $l\le N-1$ two independent strands 
$\{S^{(l)}_{k,2i-1},S^{(l)}_{k,2i},\ k=0,\dots,2^{N-l}\}$
build from $\{S^{(l-1)}_{k,i},k=0,\dots,2^{N-l+1}\}$ as in \reff{eq-not.1}. Thus, 
for each generation $l<N$, we obtain a collection of $2^l$ independent
local times $\{\{l^{(l)}_{2^{N-l},i}(x),x\in \Z^d\},\ i=1,\dots,2^l\}$,
associated with the trajectories 
$\{\{ S^{(l)}_{k,i},k=0,\dots,2^{N-l}\},\ i=1,\dots,2^l\}$.

For any $N$ and $l$, we define for $i=1,\dots,2^l$
\be{eq-not.2}
\D_i^{(N,\ l)}(z)=\{x\in \Z^d:\ l^{(l)}_{2^{N-l},i}(x)>z\},
\ee
and for $i=1,\dots,2^{l-1}$
\be{eq-not.3}
\C_i^{(N,\ l)}(z)=\{x\in \Z^d:\ \min(l^{(l)}_{2^{N-l},2i-1}(x)
,l^{(l)}_{2^{N-l},2i}(x))>z\}.
\ee
\noindent{\bf Step 1.} We first show that if $\eta=\eta'+\eta''$, then for any $\delta
\in ]0,1[$, and $l<N-1$
\be{eq-not.6}
\{\sum_{i=1}^{2^{l}} |\D_i^{(N,\ l)}(z)|>\eta\}\subset
\{\sum_{i=1}^{2^{l+1}} |\D_i^{(N,\ l+1)}\pare{(1-\delta)z}|>\eta''\}\cup
\{\sum_{i=1}^{2^{l}} |\C_i^{(N,\ l+1)}\pare{\delta z}|>\eta'\}.
\ee
We first fix one strand $\{S^{(l)}_{k,i},k=0,\dots,2^{N-l}\}$ at generation $l$.
To lighten notations, we set $m=2^{N-l-1}$. Then, on 
$\{l^{(l)}_{2m,i}(x)>z,\ S^{(l)}_{m,i}=\tilde x\}$ (for $\tilde x\in \Z^d$),
we have 
\[
z<l^{(l)}_{2m,i}(x)\le l^{(l+1)}_{m,2i-1}(\tilde x-x)+l^{(l+1)}_{m,2i}(\tilde x-x).
\]
Thus, we have either of the two following possibilities on
$\{l^{(l)}_{2m,i}(x)>z,\ S^{(l)}_{m,i}=\tilde x\}$: for $0<\delta<1$
\begin{itemize}
\item[(i)] $\max(\ l^{(l+1)}_{m,2i-1}(\tilde x-x),
\ l^{(l+1)}_{m,2i}(\tilde x-x))>(1-\delta)z$.
\item[(ii)] $\min(\ l^{(l+1)}_{m,2i-1}(\tilde x-x),
\ l^{(l+1)}_{m,2i}(\tilde x-x))>\delta z$.
\end{itemize}
Thus, by partitioning over $\{S^{(l)}_{m,i}=\tilde x,\ \tilde x\in \Z^d\}$, we obtain
\ba{eq-not.7}
\{x:l^{(l)}_{2m,i}(x)>z\}\subset \bigcup_{\tilde x\in \Z^d} 
\{S^{(l)}_{m,i}=\tilde x\}&\bigcap&
\big(\{x: \max(l^{(l+1)}_{m,2i-1}(\tilde x-x),
\ l^{(l+1)}_{m,2i}(\tilde x-x))>(1-\delta)z\}\cr
&&\cup\{x: \min(l^{(l+1)}_{m,2i-1}(\tilde x-x),\ l^{(l+1)}_{m,2i}(\tilde x-x))>\delta z\}\big)
\ea
Thus, by taking the cardinal of each set, we obtain for $i=1,\dots,2^l$,
\be{eq-not.8}
\begin{split}
&|\D_i^{(N,\ l)}(z)|\le  \sum_{\tilde x\in \Z^d} 1\{S^{(l)}_{m,i}=\tilde x\}\Big(
\sum_{j=2i-1,2i}|\{x:l^{(l+1)}_{m,j}(\tilde x-x)>(1-\delta)z\}|\\
&\qquad\qquad+ |\{x:\min(l^{(l+1)}_{m,2i-1}(\tilde x-x),
\ l^{(l+1)}_{m,2i}(\tilde x-x))>\delta z\}|\Big)\\
&\qquad\le \!\!\! \sum_{j=2i-1,2i}\!\!\!|\{x:l^{(l+1)}_{m,j}(x)>(1-\delta)z\}|+
|\{x:\min(l^{(l+1)}_{m,2i-1}(x),\ l^{(l+1)}_{m,2i}(x))>\delta z\}|\\
&\qquad\le  |\D_{2i-1}^{(N,\ l+1)}\pare{(1-\delta)z}|
+|\D_{2i}^{(N,\ l+1)}\pare{(1-\delta)z}|+
|\C_i^{(N,\ l+1)}\pare{\delta z}|.
\end{split}
\ee
Thus, \reff{eq-not.6} follows at once.

\noindent{\bf Step 2.} We show now that if we partition the size of
the level-set $\eta$ into $\eta'+\eta''$, then
\be{eq-not.11}
\acc{\sum_{i=1}^{2^{l}} |\D_{i}^{(N,\ l)}(z)|>\eta'+\eta''}\subset
\acc{\sum_{i=1}^{2^{l+1}} |\D_{i}^{(N,\ l+1)}\pare{(1-\delta)z}|>\eta''} \cup \A
\ee
with,
\[
\A:=\acc{
\sum_{i=1}^{2^{l}}l^{(l+1)}_{2^{N-l-1},2i}\pare{\D_{2i-1}^{(N,\ l+1)}(\delta z)}
\ge \delta z \eta'}.
\]
Indeed, note that $\C_i^{(N,\ l+1)}\pare{\delta z}$ are the sites of 
$\D_{2i-1}^{(N,\ l+1)}(\delta z)$ 
where $l^{(l+1)}_{2^{N-l-1},2i}\ge\delta z$. Thus,
\be{eq-not.9}
l^{(l+1)}_{2^{N-l-1},2i}\pare{\D_{2i-1}^{(N,\ l+1)}(\delta z)}\ge\delta z
|\C_i^{(N,\ l+1)}\pare{\delta z}|,
\ee
so that
\be{eq-not.10}
\sum_{i=1}^{2^{l}}l^{(l+1)}_{2^{N-l-1},2i}\pare{\D_{2i-1}^{(N,\ l+1)}(\delta z)}>
\delta z \sum_{i=1}^{2^{l}} |\C_i^{(N,\ l+1)}\pare{\delta z}|,
\ee
and we deduce Step 2 from \reff{eq-not.6} and \reff{eq-not.10}.

\noindent{\bf Step 3} We partition further \reff{eq-not.11} to get rid
of the event that one of the $\D_{2i-1}^{(N,\ l+1)}(\delta z)$ in $\A$
be too large. Thus, for an arbitrary positive constant $a$ to be chosen later,
\be{eq-not.12}
\A \subset\pare{\A\quad\bigcap_{i=1}^{2^l} \acc{|
\D_{2i-1}^{(N,\ l+1)}(\delta z)|\le \frac{\eta}{\delta^a}}}
\bigcup_{i=1}^{2^l} \{| \D_{2i-1}^{(N,\ l+1)}(\delta z)|>\frac{\eta}{\delta^a}\}.
\ee
Now, we denote
\be{eq-not.4}
A^{N}_{l}(z,\eta)=P\pare{\sum_{i=1}^{2^l} |\D_i^{(N,\ l)}(z)|>\eta},
\ee
and,
\be{eq-not.5}
B^{N}_{l}(\eta;(z,w))=P\pare{\sum_{i=1}^{2^{l-1}} l^{(l)}_{2^{N-l},2i}
\pare{\D_{2i-1}^{(N,\ l)}(z)}>\eta;\quad \forall i=1,
\dots,2^{l-1},\quad|\D_{2i-1}^{(N,\ l)}(z)|<w}.
\ee
We consider a decomposition of $\eta$ into $N-1$ positive numbers 
$\eta_1,\dots,\eta_{N-1}$,
and we denote $\bar \eta_i=\eta_{i+1}+\dots+\eta_{N-1}$. 
Now, at generation $l<N-1$, we apply Step 1 and Step 2 and
\reff{eq-not.12} with $\eta=\bar\eta_l$, $\eta'=\eta_l$
and $\eta''=\bar \eta_{l+1}$. If we further
take averages on both sides of \reff{eq-not.12}, we obtain
\be{eq-not.13}
A^{N}_{l}(z,\bar \eta_{l})\le 
A^{N}_{l+1}((1-\delta)z,\bar \eta_{l+1})+2^l A^{N-l-1}_{0}(\delta z,
\frac{\eta}{\delta^a})+ B^{N}_{l+1}(\delta z \eta_{l+1};
(\delta z,\frac{\eta}{\delta^a})).
\ee

Now, we define $\Theta(z,\eta)=(\delta z, 
\frac{\eta}{\delta^a})$, $\Gamma(z,\eta)=((1-\delta)z,\eta)$ and
for each $l\le N$, $m_l(z,\eta)=(\frac{\delta}{l} z \eta;\Theta(z,\eta))$. 
By iterating
\reff{eq-not.13} for the term $A^{N}_{l+1}$ (until $l=N-1$ since $A^k_k=0$),
and choosing $\eta_i=\eta/(N-1)$ for $i=1,\dots,N-1$, we obtain
\be{eq-not.15}
A^{N}_0(z,\eta)\le \sum_{l=1}^{N-1}\acc{2^{l-1}A^{N-l}_0\circ\Theta+
B^N_l\circ m_{N-1}}\circ \Gamma^{l-1}(z,\eta) .
\ee
On the right hand side of \reff{eq-not.15}, we have desirable $B$-terms,
and $A^k_0$-terms which we get rid off by iterating \reff{eq-not.15}. Note that
the action of iterates of $\Gamma$ on $(z,\eta)$
will be innocuous as we choose later $\delta$ very small; however, the action of
$\Theta$ must be traced carefully. Thus, in \reff{eq-not.13}, we say that in the $A$-terms
of the right hand side, $\Theta$ acts once. Also, a given $A$-term, say $A_0^{N-l}$ has
argument $\Theta\circ\Gamma^{l-1}(z,\eta)=(\delta(1-\delta)^{l-1}z,\frac{\eta}{\delta^a})$,
and in the induction, we need to decompose $\frac{\eta}{\delta^a}$ into $N-l-1$ equal parts so
as to obtain
\[
A^{N-l}_0\Theta\circ\Gamma^{l-1}(z,\eta)\le \sum_{l'=1}^{N-l-1}\acc{2^{l'-1}A^{N-l-l'}_0
\circ\Theta^2+ B^{N-l}_{l'}\circ m_{N-l-1}\circ\Theta}\circ \Gamma^{l-1+l'-1}(z,\eta) .
\]
We describe now in more details the $B$-terms we eventually obtain. In a {\it generic}
$B$-term, let $\nu\ge 0$ be the number of times $\Theta$ has acted, and for $i=1,\dots,\nu$,
let $l_i$ be the number of times $\Gamma$ has acted between the $(i-1)^{th}$ and
$i^{th}$ action of $\Theta$, and let $l\ge 1$ be the number of times $\Gamma$ acts
after the $\nu$-actions of $\Theta$. We assume $1\le l_1+\dots+l_{\nu}+l\le N$.
We set 
\[
k=l_1+\dots+l_{\nu},\quad\text{and}\quad k''=k-\nu+l-1. 
\]
For a single
choice $(\nu,\ l_1,\dots,\ l_{\nu},\ l)$, we have $2^{k-\nu}$ $B$-terms of the form
$B^{N-k}_l$ and with argument $m_{N-k}\circ\Theta^{\nu}\circ \Gamma^{k''}(z,\eta)$.
Note that the total number of $B$-terms labelled $B^{N-k}_l$ is the same as those
labelled $B^{N-k}_0$, a number we call $c(k)$ which is easily seen from
\reff{eq-not.15} to satisfy
\be{eq-not.16}
c(k)\le 2^0 c(k-1)+2^1 c(k-2)+\dots+2^{k-1} c(0),
\quad\text{with}\quad
c(0)=1,\ c(1)=1,\ c(2)=3,\quad\text{etc}...
\ee
Now, an immediate induction shows that \reff{eq-not.16} imposes the bound
$c(k)\le 2^{2k}$. Thus, we obtain
\be{eq-not.17}
A^{N}_0\le \sum_{k,l:\ N\ge k+l}\!\!\! 2^{2k} \sup_{k''\le k+l}\quad
\sup_{\nu< k} \quad B^{N-k}_l\circ m_{N-k}\circ\Theta^{\nu}\circ \Gamma^{k''}
\ee
We write in details the $B$-term in \reff{eq-not.17} for
a choice of $(\nu',\ l_1,\dots,\ l_{\nu},\ l)$, and $\nu=\nu'+1$.
\be{eq-not.20}
B^{N-k}_l\circ m_{N-k}\circ\Theta^{\nu'}\circ \Gamma^{k''}(z,\eta)=
P\pare{ \sum_{i=1}^{2^{l-1}} l^{(l)}_{2^{N-k-l},2i}\pare{D_{2i-1}}>
\frac{\delta(1-\delta)^{k''}}{N-k} \frac{\delta^\nu z\eta}{\delta^{\nu a}},\ 
\bigcap_{i=1}^{2^{l-1}}\G_i },
\ee
where, for $i=1,\dots, 2^{l-1}$, we used the shorthand notations
\[
D_{2i-1}:=\D_{2i-1}^{(N-k,\ l)}\pare{ (1-\delta)^{k''} \delta^{\nu} z},\quad
\text{and}\quad
\G_i:= \acc{ |D_{2i-1}|< \frac{\eta}{\delta^{\nu a}} }.
\]

We take now $a=d/(d-2)$. To understand this choice, 
note that we deal in \reff{eq-not.20} with a sum
of $2^{l-1}$ independent terms whose tail distribution is controlled by
inequality \reff{lem-ACa}. It will turn out, for the forthcoming choice 
of $(z,\eta)$, that the sum in \reff{eq-not.20} behaves similarly as one of its term. 
Now, if we were asking for the probability that
\[
l^{(l)}_{2^{N-k-l},2i}\pare{D_{2i-1}}>\frac{\delta^\nu z\eta}{\delta^{\nu a}},
\quad\text{with}\quad \left|D_{2i-1}\right|<\frac{\eta}{\delta^{\nu a}},
\]
then, estimates \reff{lem-ACa} would give a bound $\exp(-\kappa_d
F\pare{ \frac{ \delta^{\nu} z\eta}{\delta^{a\nu}}, \frac{\eta}{\delta^{a\nu}} })$. 
Thus, $\Theta^{\nu}$ will not ruin the use of estimate \reff{lem-ACa} if
for the function $F$ given in \reff{lem-ACa} we have that
$F\pare{ \frac{ \delta^{\nu} z\eta}{\delta^{a\nu}}, \frac{\eta}{\delta^{a\nu}} }$
is independent of $\delta$. This is what we achieve by choosing $a=d/(d-2)$.

\noindent{\bf Step 4.} 
We are now ready to evaluate the level sets distribution. Note that
\be{level-D}
\acc{ |\LL_j|>\frac{2^N y_j}{N^{2\alpha_{j+1}}} }\subset
\acc{ |\D^{(N,0)}_1(N^{\alpha_j})|>\frac{2^N y_j}{N^{2\alpha_{j+1}}} }.
\ee
We rewrite the $B$-term of \reff{eq-not.20} with
$z=N^{\alpha_j}$ and $\eta=2^N y_j/N^{2\alpha_{j+1}}$.
\be{key.1}
B^{N-k}_l(.)=P\pare{ \sum_{i=1}^{2^{l-1}} X_i^{(l)}>x_N,\quad \bigcap_{i=1}^{2^{l-1}}
 \G_i},
\ee
with 
\be{key.2}
X_i^{(l)}=\delta^{\nu(a-1)}\quad l^{(l)}_{2^{N-k-l},2i}\pare{D_{2i-1}},
\quad\text{with}\quad
D_{2i-1}:=\D_{2i-1}^{(N-k,\ l)}\pare{ (1-\delta)^{k''} \delta^{\nu} N^{\alpha_j}},
\ee
and, as we chose $y_j=y/(2M_N)$ for $j>0$,
\be{key.3}
x_N=\frac{\delta(1-\delta)^{k''}}{N-k} \frac{N^{\alpha_j}}{N^{2\alpha_{j+1}}2M_N}2^Ny ,
\quad\text{and}\quad
\G_i=\{|D_{2i-1}|<\frac{ 2^N y}{ N^{2\alpha_{j+1}}2M_N \delta^{\nu a}}\}.
\ee
For $B^{N-k}_l$ to be small, we need $2^{l-1}E[X_i^{(l)}1\{\G_i\}]<x_N/2$.
Thus, we show in Lemma~\ref{lem-app2}, that there is a constant $C_0$ such that
\be{key.4}
E[X_i^{(l)}]\le C_0\ \delta^{\nu(a-1)}
\ 2^{\frac{2}{3}(N-k-l)}\exp\pare{-\kappa_3 \frac{2}{3}
\delta^{\nu}(1-\delta)^{k''} N^{\alpha_j} }.
\ee
Now, to get rid of the term $(1-\delta)^{k''}$ we take $\delta=1/N$ (since
$k''\le 2N$). Now, recall that when $d=3$, then $a=3$. Thus,
we have $2^{l-1}E[X_i^{(l)}]<x_N/2$, if for some constant $C_1$ 
\be{key.5}
2^l 2^{\frac{2}{3}(N-k-l)} \delta^{2\nu}\ \exp\pare{-c_0 \delta^{\nu}N^{\alpha_j} }\le
\frac{C_1\delta}{N} \frac{2^N N^{\alpha_j} y}{N^{2\alpha_{j+1}}M_N },
\ee
where we set $c_0:=\sup_N\acc{\frac{2}{3}\kappa_3(1-\frac{1}{N})^{2N}}>0$.
If we set $x= \delta^{\nu} N^{\alpha_j} $, then \reff{key.5} holds as soon as
\be{key.6}
x^2 \exp\pare{ -c_0 x}\le \frac{C_1 y}{N^3} \frac{N^{\alpha_j}}
{N^{2(\alpha_{j+1}-\alpha_j)}} 2^{(N-l)\frac{1}{3}}.
\ee
Since $\alpha_{j+1}-\alpha_j=\alpha$ can be chosen arbitrarily small,
\reff{key.6} follows as soon as $\uchi>3$. Thus, if we set $Y_i^{(l)}=
X_i^{(l)}1\{\G_i\}$ and $\bar Y_i^{(l)}=Y_i^{(l)}-E[Y_i^{(l)}]$, 
we have
\be{key.7}
P\pare{ \sum_{i=1}^{2^{l-1}} X_i^{(l)}> x_N,\quad \bigcap_{i=1}^{2^{l-1}} \G_i}\le
P\pare{ \sum_{i=1}^{2^{l-1}} \bar Y_i^{(l)}> \frac{x_N}{2}}.
\ee
We have reached now a large deviation estimate for which
Lemma~\ref{App-UB} is devised. We first need tail estimates for $Y_i^{(l)}$.

\noindent{{\bf Step 5}}: To obtain tail estimates, we rely on
Lemma 1.2 of \cite{ACa}, 
\ba{key.8}
P\pare{Y_i^{(l)}>u}&=&
P\pare{ X_i^{(l)}>u,\ \G_i}\le E\cro{ \exp\pare{ -\frac{u \delta^{-\nu(a-1)}}
{|\D_{2i-1}|^{2/d}} } 1\acc{\G_i} }\cr
&\le & \exp\pare{ -\xi_N u},\quad\text{with}\quad
\xi_N=\pare{ \frac{N^{2\alpha_{j+1}}M_N}{2^N y}}^{2/d}.
\ea
We show now that for $\uchi\ge 2$ (and $\alpha<1/2$), we have $2^{4/3l}\xi_N^2
E[(Y_i^{(l)})^2]\le 1$. Using Lemma~\ref{lem-app2}, there is a constant $c$
\ba{key.13}
2^{4/3l}\xi_N^2E\cro{(X_i^{(l)})^2}&=&
2^{4/3l}\pare{ \frac{N^{2\alpha_{j+1}}M_N}{2^N y}}^{4/d} \delta^{2\nu(a-1)}
2^{\frac{4}{3}(N-k-l)} \exp(-c \delta^{\nu} N^{\alpha_{j}})\cr
&\le &
\pare{ \frac{N^{2(\alpha_{j+1}-\alpha_{j})}M_N}{y}}^{4/3}
\frac{1}{N^{(4-8/3)\alpha_j}} (\delta^{\nu} N^{\alpha_{j}})^4
\exp(-c \delta^{\nu} N^{\alpha_{j}}).
\ea
The right hand side of \reff{key.13} can be made smaller than 1
if $(2\alpha+1)4/3<(4-8/3)\alpha_j$, i.e. as $\uchi\ge 2$ and $\alpha<1/2$.
Thus, Lemma~\ref{App-UB} with the choice $\gamma=1/4$ yields
\be{key.11}
P\pare{ \sum_{i=1}^{2^{l-1}} \bar Y_i^{(l)}>\frac{x_N}{2}}\le
\exp\pare{ c_u\pare{2^{4/3l}\xi_N^2 \gamma^2E\cro{(Y_i^{(l)})^2} }^{1-\gamma}-
\frac{\gamma\xi_N x_N}{4} }.
\ee
Thus, we obtain that for some constant $c>0$ and $N$ large,
\ba{key.16}
P\pare{ \sum_{i=1}^{2^{l-1}} \bar X_i^{(l)}>
\frac{x_N}{2},\quad \bigcap_{i=1}^{2^{l-1}}\G_i}
&\le& c\exp\pare{- \frac{\xi_N x_N}{16}}\le c\exp\pare{ -(2^Ny)^{1/3} N^{\zeta} },
\ea
with 
\be{key.17}
\zeta=\frac{\alpha_{j}}{3}-2(\alpha_{j+1}-\alpha_{j})-2-\frac{1}{3}.
\ee
Thus, $\zeta>0$ as soon as $\uchi>7$, and $\alpha$ small enough.

\subsection{Proof of Upper Bound in \reff{UB-prop1}}\label{tilt}
Note that in dimension 3, we are left with showing that for
$\LL_0:=\{x:\ 0<l_{2^N}(x)<N^{\uchi}\}$, we have for $\bar c>0$ and $y_0=y/2$
\[
P\pare{\sum_{\LL_0}l_{2^N}^2(x)>y_0}\le \exp\pare{-\bar c y^{1/3}n^{1/3}}.
\]
The approach is close to the proof of Lemma 3.1 in \cite{ACb}.
However, in order to get rid of a logarithmic term, inherent in
the proof in \cite{ACb}, additional work is needed. On the other hand, the
proof we now present does not work in dimensions $d\ge 4$. 

We keep the notations of Section~\ref{sec-new}.
\be{decomp.aa1}
\sum_{x \in \LL_0} l_n^2(x)\leq n+1 + 2 Z^{(0)},\quad\text{with}\quad
Z^{(0)}=\sum_{x \in \LL_0,} \sum_{0 \leq k < k' \leq 2^N} 
1\{S_k^{(0)}=S_{k'}^{(0)}=x\} \, .
\ee
Now,
\ba{decomp.aa2}
Z^{(0)} & \leq &
\sum_{x}1\{l_{2^{N-1}}^{(0)}(x) \leq N^{\uchi}\} \, \sum_{0 \leq k < k' \leq 2^{N-1}}
1\{S_k^{(0)}=S_{k'}^{(0)}=x\}\cr
&& + \sum_{x}1\{l_{2^N}^{(0)}(x)-l_{2^{N-1}}^{(0)}(x) \leq N^{\uchi}\} \,
\sum_{2^{N-1} \leq k < k' \leq 2^{N}} 1\{S_k^{(0)}=S_{k'}^{(0)}=x\}\cr
&& + \sum_{x}1\{l_{2^{N-1}}^{(0)}(x) \leq N^{\uchi}\} \,
\sum_{0 \leq k \leq 2^{N-1} \leq  k' \leq 2^{N}} 1\{S_k^{(0)}=S_{k'}^{(0)}=x\}\cr
&\le&Z^{(1)}_1 + Z^{(1)}_2 + J^{(1)}_1 \, .
\ea
where we have defined for $i=1$ and $i=2$
\[
Z^{(1)}_i = \sum_{x}1\{ l_{2^{N-1},i}^{(1)}(x) \leq N^{\uchi}\} \, 
		\sum_{0 \leq k < k' \leq 2^{N-1}}
		1\{S_{k,i}^{(1)}=S_{k',i}^{(1)}=x\} \, ,
\]
and the intersection times of the two independent strands over 
$\{ l_{2^{N-1},1}^{(1)}(x) \leq N^{\uchi}\}$ is
\[
J^{(1)}_1 = \sum_{x}1\{ l_{2^{N-1},1}^{(1)}(x) \leq N^{\uchi}\}
		l_{2^{N-1},1}^{(1)}(x) l_{2^{N-1},2}^{(1)}(x) \, .
\]
Iterating this procedure, we get 
\be{legall-dec}
Z^{(0)} \leq \sum_{l=1}^{N-1} \sum_{k=1}^{2^{l-1}} J^{(l)}_k \, ,
\ee
where for each $l \in \{1,\cdots, N-1\}$, the random variables 
$\{J^{(l)}_k; 1 \leq k \leq 2^{l-1}\}$ are i.i.d.\ , with
\[
J_k^{(l)}=\sum_{x}1\{l_{2^{N-l},2k-1}^{(l)}(x) \leq N^{\uchi} \} 
l_{2^{N-l},2k-1}^{(l)}(x) l_{2^{N-l},2k}^{(l)}(x),
\]
We now introduce a partition of $\{l_{2^{N-l},2k-1}^{(l)}(x) \leq N^{\uchi} \}$
in terms of
\be{turb-5}
\D^{(l)}_{k,i}=\{x:\  N^{\chi_i}\le l^{(l)}_{2^{N-l},2k-1}(x)< N^{\chi_{i+1}}\},
\quad\text{with}\quad \chi_i=\frac{\chi(1+\delta)^i}{\log(N)},\quad
i=0,\dots,M_N.
\ee
We choose $\chi$ such that $N^{\chi_0}=1$, and $\delta<1/3$. The reason for 
such choices will become clear later. Note that $M_N$ is of order
$\log(\log(N))$. Also, we introduce for $k=1,\dots,2^{l-1}$
\be{def-Jik}
J^{(l)}_{k,i} = \sum_{x}1\{ x\in \D^{(l)}_{k,i}\}
                l^{(l)}_{2^{N-l},2k-1}(x) l^{(l)}_{2^{N-l},2k}(x),\quad\text{and}\quad
J^{(l)}_{k}=\sum_{i=0}^{M_N} J^{(l)}_{k,i}.
\ee
Finally, we need the self-intersections of the $2^l$ strands at generation $l$
\be{def-Zlk}
Z^{(l)}_k=\sum_{x }
 \sum_{0 \leq m < m'\leq 2^{N-l}} \!\!\! 1\{ l^{(l)}_{2^{N-l},k}(x) \leq N^{\uchi}\}
1\acc{S^{(l)}_{m,k}=S^{(l)}_{m',k}=x},
\ee
We bootstrap a little differently than in the proof of Lemma 2.1 of \cite{ACb}.
Thus, at each generation $l$, and for level-set index $i$, we introduce
the good-sets 
\[
\forall k=1,\dots,2^{l-1},\quad\forall i=0,\dots,M_N,\qquad
\G^{(l)}_{k,i}=\{ |\D^{(l)}_{k,i}|<\frac{4 y_0 2^N}{N^{2\chi_i}}\},\quad\text{and}\quad
\G^{(l)}=\bigcap_{k,i} \G^{(l)}_{k,i}.
\]
As in equation (35) of \cite{ACb}, we have
\be{good-set}
(\G^{(l)})^c \subset
\{Z^{(l)}_1+\dots + Z^{(l)}_{2^l}> y_0 2^N\}.
\ee
It is important to note that contrary to (35) of \cite{ACb}, we have kept
the threshold $y_0 2^N$. Thus,
\ba{decompZ0}
P\pare{Z^{(0)} > y_0 2^N}&\le & 
P\pare{\sum_{l=1}^{N-1} \sum_{k=1}^{2^{l-1}} J^{(l)}_k >y_0 2^N}\cr
&\le & P\pare{\sum_{l,k} J^{(l)}_k>y_0 2^N, \bigcap_{l=1}^{N-1} \G^{(l)} }+
\sum_{l=1}^N P\pare{ (\G^{(l)})^c }\cr
&\le & P\pare{\sum_{l,k} J^{(l)}_k>y_0 2^N, \bigcap_{l=1}^{N-1} \G^{(l)} }+
\sum_{l=1}^N P\pare{ \sum_{j=1}^{2^l} Z^{(l)}_j> y_0 2^N}\, .
\ea
Now, by writing self-intersection in terms of intersection of independent
strands, and proceeding by induction, we obtain
\ba{no-weak}
P\pare{Z^{(0)} > y_0 2^N}\!\!\!&\le&\!\!\!\sum_{L=1}^N 2^{L-1}
P\!\!
\pare{\sum_{l=L}^{N-1} \sum_{k=1}^{2^{l-1}} J^{(l)}_k >y_0 2^N, 
\bigcap_{l=L}^N \G^{(l)}\!\!}
\!\!\!+2^{N-1}P\!\!\pare{ \sum_{k=1}^{2^{N-1}} Z^{(N-1)}_k>y_0 2^N \!\!}
\cr
\!\!\!&\le &\!\!\! \sum_{L=1}^N 2^{L-1}
P\!\!\pare{\sum_{l=L}^{N-1}\sum_{i=0}^{M_N}
\sum_{k=1}^{2^{l-1}} J^{(l)}_{k,i}1\{\G^{(l)}_{k,i}\} >y_0 2^N}.
\ea
The last term of the first line in \reff{no-weak} has vanished since
$Z^{(N-1)}_k\le 1$, and we choose $y_0> 1/2$. 
Also, note that \reff{no-weak} is different from inequality (36) of \cite{ACb} in
having the sum over $l$ inside the probability.
Now, Lemma~\ref{lem-app2} of the Appendix allow us to center the $J^{(l)}_k$,
since
\[
E \cro{\sum_{l=1}^{N} \sum_{k=1}^{2^{l-1}} J^{(l)}_k} 
\leq C_3 2^N. 
\]
Actually, we rather need to center $Y^{(l)}_{k,i}:=J^{(l)}_{k,i}1\{\G^{(l)}_{k,i}\}$.
Thus, let $\bar Y^{(l)}_{k,i}=Y^{(l)}_{k,i}-E[Y^{(l)}_{k,i}]$,
and if we set $\bar y<(y_0-1)/2-C_3$, and choose $y$ large enough so that 
$\bar y>y_0/2$,
\[
P \pare{\sum_{x \in \LL_0} l_n^2(x) \geq 2^Ny_0}
\leq \sum_{L=1}^N 2^{L-1}
P \pare{ \sum_{l=L}^{N-1}
\sum_{i=0}^{M_N}  \sum_{k=1}^{2^{l-1}} \bar Y^{(l)}_{k,i} \geq  \bar y 2^N}
\, .
\]
Now, fix $L$ and note that for any sequences $\{q_i,p_l^{(i)},i=0,\dots,M_N,\ 
l=L,\dots,N\}$ with
\[
\sum_i q_i\le 1,
\quad\text{and}\quad \sum_l p_l^{(i)}\le 1, \quad\text{for any }i=0,\dots,M_N,
\]
we have
\be{no-weak1}
P \pare{ \sum_{l=L}^{N-1}\sum_{i=0}^{M_N} 
  \sum_{k=1}^{2^{l-1}} \bar Y^{(l)}_{k,i} \geq  \bar y 2^N}
\le \sum_{l=L}^{N-1} \sum_{i=0}^{M_N} P\pare{\sum_{k=1}^{2^{l-1}} \bar Y^{(l)}_{k,i} 
 \geq  p_l^{(i)}q_i\bar y 2^N}
\ee
In order to use Lemma~\ref{App-UB}, we 
need exponential estimates for the $\bar Y_{k,i}^{(l)}$. 
Note first that
\be{Jki}
J_{k,i}^{(l)}\le N^{\chi_{i+1}} l^{(l)}_{2^{N-l},2k}\pare{ \D^{(l)}_{k,i} },
\ee
We use Lemma 1.2 of \cite{ACa} to obtain
\[
P\pare{  Y_{k,i}^{(l)}>u}=
P\pare{ J_{k,i}^{(l)}>u, \G_{k,i}^{(l)} }\le c_3
E\cro{ \exp\pare{-
\kappa_3\frac{u}{N^{\chi_{i+1}}|\D^{(l)}_{k,i}|^{2/3} }}1\{ \G_{k,i}^{(l)}\} }.
\]
We have two bounds on $|\D^{(l)}_{k,i}|$: either we recall that, on $\G_{k,i}^{(l)}$,
the volume is bounded by $2y 2^N/N^{2\chi_i}$, or the trivial bound by the total
time $2^{N-l}$. Thus,
\be{no-weak2}
P\pare{  Y_{k,i}^{(l)}>u}\le C_3 \exp\pare{ -\xi^{(l)}_{i} u}
\quad\text{with any}\quad
\xi^{(l)}_{i}\le\frac{\kappa_3}{N^{\chi_{i+1}}2^{\frac{2}{3}N}}
\max\pare{2^l,\pare{ \frac{N^{2\chi_i}}{2y} } }^{2/3}.
\ee
We define $\xi_N=1/2^{\frac{2}{3}N}$, and for a fixed $i$, we choose for convenience
(for a $\delta<1/3$)
\be{no-weak3}
\frac{\xi^{(l)}_{i}}{\xi_N}=\frac{\kappa_3}{N^{\chi_{i+1}}}
\left\{ \begin{array}{ll}
                \frac{N^{\frac{4}{3}\chi_i}}{(2y)^{2/3} }&
\mbox{ for } l\le l^*_i \, ,
                \\
                2^{\frac{1}{6}l} & \mbox{ for } l> l^*_i \, ,
                \end{array}
        \right.
\ee
where $l^*_i$ is such that 
\[
2^{\frac{1}{6}l^*_i}=\frac{N^{\frac{4}{3} \chi_i} }{(2y)^{2/3}},\quad\text{so that}\quad
l^*_i=\pare{\frac{8}{\log(2)} \chi (1+\delta)^i-\frac{2}{3}\log(2y)}^+.
\]
We wish now to use Lemma~\ref{App-UB}, or rather Remark~\reff{app-rem1}, with
$\Gamma=\xi^{(l)}_{i}$ and $X_i=Y_{k,i}^{(l)}$. Thus, we first bound $\Gamma^2
E[X_i^2]$ using Lemma~\ref{lem-app3} and \reff{Jki}
\ba{no-weak10}
(\xi^{(l)}_{i})^2 E[(Y_{k,i}^{(l)})^2]&\le & 
(\xi^{(l)}_{i})^2 E[(J_{k,i}^{(l)})^2]\le C\frac{ 2^{\frac{2}{6}l}
N^{\frac{8}{3} \chi_i} }{(2y)^2 2^{\frac{4}{3} N} }2^{\frac{4}{3}(N-l)} e^{
-\kappa_3 N^{\chi_i}}\cr
&\le & C\frac{ \pare{ N^{\chi_i}}^{8/3} e^{-\kappa_3 N^{\chi_i}}}{ (2y)^2}\frac{1}{2^l}.
\ea
For $\gamma>0$ small, we denote for convenience $\sigma_N=\gamma^2\Gamma^2E[X_i^2]$.
By Lemma~\ref{App-UB}, we obtain for any $\gamma\in ]0,1[$
\be{no-weak11}
P\pare{ \sum_{k=1}^{2^{l-1}} \bar Y_{k,i}^{(l)}\ge p_l^{(i)}q_i y 2^N }
\le
\exp\pare{ -\frac{\gamma \xi^{(l)}_{i}}{2}p_l^{(i)}q_i \bar y 2^N +c_u 2^l
\max\pare{\sigma_N,\sigma_N^{1-\gamma}}}.
\ee
Assume now that we can choose $p_l^{(i)}$ and $q_i$ such that for some constant $c$,
and $y$ large (but fixed as $N$ tends to infinity)
\be{norm}
\frac{\xi^{(l)}_{i}}{\xi_N} p_l^{(i)}q_i \ge \frac{c}{y^{2/3}}.
\ee
Then, \reff{no-weak11} yields the upper bound in \reff{UB-prop1} if 
\be{no-weak13}
c\gamma y^{1/3} 2^{\frac{1}{3}N}\ge 8 c_u 2^l \max\pare{\sigma_N,\sigma_N^{1-\gamma}}. 
\ee
Note that from \reff{no-weak10}
\[
\sigma_N\le \frac{C}{(2y)^2}\sup_{x\ge 0} \pare{ x^{8/3} e^{-\kappa_3 x}}\frac{1}{2^l},
\]
so that \reff{no-weak13} holds if $2^{\frac{1}{3}N}\gg 2^{\gamma l} $.
Now, since $l\le N$,
\reff{no-weak13} holds as soon as $\gamma<1/3$ for $N$ large enough.

Finally, we choose $p_l^{(i)}$ and $q_i$ to fulfill \reff{norm}.
We set $\alpha:=\frac{1}{2}(\frac{1}{3}-\delta) \chi$, and
\be{no-weak15}
q_i=q\pare{ \frac{N^{\chi_{i+1}}}{N^{\frac{4}{3}\chi_i}} }^{1/2}=
q\exp\pare{ -\alpha(1+\delta)^i},\quad\text{with $q$ such that }\quad
\sum_{i=0}^{M_N} q_i\le 1.
\ee
Note that it is possible to find such a $q$ which depends on $\chi$ and $\delta$.
Now, fix $i$, and choose
\be{no-weak16}
\forall l\le l^*_i,\quad p_l^{(i)}=p_i^* \exp\pare{ -\alpha(1+\delta)^i},\quad
\text{whereas if}\quad l>l^*_i,\quad p_l^{(i)}=\bar p_i 
\frac{N^{\chi_{i+1}} e^{\alpha(1+\delta)^i}}{2^{\frac{1}{6}l}},
\ee
with two normalizing constants $p_i^*$ and $\bar p_i$ to be chosen later.
Note that for $l>l^*_i$ 
\be{bound-pi}
\frac{N^{\chi_{i+1}} e^{\alpha(1+\delta)^i} }{2^{\frac{1}{6}l}}\le
\pare{ \frac{N^{\frac{4}{3}\chi_{i}}N^{\chi_{i+1}}}
{2^{\frac{1}{6}l_i^*}2^{\frac{1}{6}l}} }^{1/2}\le 
\pare{\frac{(2y)^{2/3}N^{(1+\delta)\chi_{i}}}{2^{\frac{1}{6}l}} }^{1/2}.
\ee
Note that from the definition of $l_i^*$, and the choice $\delta<1/3$, we have for $l>l^*_i$
\be{bound-pi2}
N^{(1+\delta)\chi_{i}}\le N^{\frac{4}{3}\chi_{i}}=(2y)^{2/3}2^{\frac{1}{6}l_i^*}
\quad\Longrightarrow\quad
p_l^{(i)}\le \bar p_i (2y)^{2/3}  2^{-(l-l_i^*)/12} .
\ee
To see that it is possible to choose $p_i^*$ and $\bar p_i$ such that for each $i$,
$\sum_l p_l^{(i)}\le 1$, note that
\ba{no-weak17}
\sum_{l} p_l^{(i)}&\le&p_i^* l^*_i \exp\pare{ -\alpha(1+\delta)^i}+
\bar p_i (2y)^{2/3}\sum_{l>l_i^*}  2^{-(l-l_i^*)/12}\cr
&\le & p_i^* \pare{ \frac{8\chi}{\log(2)} (1+\delta)^ie^{-\alpha(1+\delta)^i}}+
\bar p_i (2y)^{2/3}\sum_{l>0} \frac{1}{2^{l/12}}\cr
&\le & p_i^*\pare{ \frac{8\chi}{\log(2)} \sup_{x>0}\acc{x e^{-\alpha x}}}+
\bar p_i\frac{(2y)^{2/3}}{2^{1/12}-1}.
\ea
It suffices now to choose $p_i^*$ as a small constant (depending only
on $\chi$), and $\bar p_i$ as a small constant times $1/(2y)^{2/3}$.
It is easy now to check that \reff{norm} holds.
\br{rem-d4}
When dimension $d=4$, the proof of Lemma 3.1 of \cite{ACb}, with
Remark~\ref{rem-app2} to obtain centering of the $J_k^{(l)}$ variables, 
can be used to obtain the upper bound \reff{UB-prop1}.
Indeed, in \cite{ACb} dimension $d\ge 5$ was used to obtain that the first
two moments of the intersection times of two independent walks were finite.
This is actually much too strong, and a close inspection of the proof of
Lemma 3.1 of \cite{ACb} shows us that we actually only need \reff{app.6}.
We omit to repeat the proof since it is similar.
\er
\subsection{Proof of the Lower Bound in \reff{UB-prop1}}\label{lower}
The proof proceed as in (66) of \cite{ACb}, by using the comparison 
$\Sigma_n^2\ge n^2/|\RR_n|$ where we denoted by $\RR_n$ the range of the walk.
Since it is a few lines, we reproduced it for the ease of reading.
Indeed, $\Sigma_n^2\ge n^2/|\RR_n|$  follows by Jensen's inequality
\be{LB-eq1}
\pare{\frac{1}{|\RR_n|}\sum_{\RR_n} l_n(x)}^2\le
\frac{1}{|\RR_n|} \sum_{\RR_n} l_n(x)^2.
\ee
Now, if $\sigma(r)$ is the first time the walk exits a ball $B(r)$, we have
\be{LB-eq2}
\acc{\sigma(r)>n}\subset \acc{|\RR_n|< |B(r)|}\subset\acc{\Sigma_n^2>\frac{n^2}
{|B(r)|}}.
\ee
Thus, if we choose a radius $r_n$ such that
$|B(r_n)|=n/y$, then $\acc{\sigma(r_n)>n}\subset \acc{\Sigma_n^2>yn}$.
We recall now the classical estimate $\P_0(\sigma(r_n)\ge n)\ge \exp(-Cn/r_n^2)$, for
some constant $C$, and this yields the lower bound in\reff{UB-prop1}.
\section{Application of Section~\ref{sec-main} to lower bounds.}\label{sec-low}
\subsection{Proof of Proposition~\ref{prop-range}}\label{sec-range}
We assume, for simplicity, that we can divide $[0,n]$ into $k_n$ periods 
of length $|B(r_n)|$. Let $T_i=(i-1)|B(r_n)|$, and
$\RR_i:=\{0,S_{T_i+1}-S_{T_i},\dots,S_{T_{i+1}}-S_{T_i}\}$ for $i=1,\dots,
k_n$. Note that $\{\RR_i,i=1,\dots,k_n\}$ 
are independent, and that for $\epsilon_0$ small, inequality~\reff{boltup} yields
\be{range-ineq1}
P(|\RR_i|<2\epsilon_0 |B(r_n)|)\le \exp\pare{-\frac{\kappa }{(2\epsilon_0)^{1/3}} 
|B(r_n)|^{1/3} }.
\ee
Now, we introduce independent Bernoulli variables 
$X_i=1\{|\RR_i|<2\epsilon_0 |B(r_n)|\}$ for $i=1,\dots,k_n$.
We rewrite \reff{range-ineq1} with a rate $I(\epsilon_0)$ 
large when $\epsilon_0$ is small, such that 
\[
E[X_i]\le \exp(-I(\epsilon_0) |B(r_n)|^{1/3}). 
\]
By Chebychev's inequality, there is a constant $c$ 
depending on $(\delta_0/\epsilon_0)$, such that when $\delta_0<\epsilon_0$
and large $n$,
\be{range-cheb}
P\pare{ \frac{1}{k_n}\sum_{i=1}^{k_n} X_i> 1-\frac{\delta_0}{\epsilon_0} }\le
\exp\pare{ -c I(\epsilon_0)|B(r_n)|^{1/3} k_n}.
\ee
On the complementary event $\{\sum_i X_i\le(1-\delta_0/\epsilon_0)k_n\}$, 
and there
are $\frac{\delta_0}{\epsilon_0}k_n$ periods, say the {\it good} periods,
where $\{|\RR_i|\ge 2\epsilon_0 |B(r_n)|\}$. We show now that if there are
enough good periods, then a fraction of the sites of $B(r_n)$ are visited
a fraction of the time $n/|B(r_n)|$. In other words,
\ba{range-idea}
\{\sigma(r_n)>n\} \cap \{|\{i\le k_n:&&\!\!\!|\RR_i|\ge 2\epsilon_0 |B(r_n)|\}|>
\frac{\delta_0}{\epsilon_0} k_n \}\cr
\subset&&\!\!\! \acc{
|\acc{x:\ l_n(x)\ge \delta_0 k_n}|\ge \epsilon_0|B(r_n)|}.
\ea
We take an issue in the left hand event in \reff{range-idea}, and by way
of contradiction, we assume that more than $(1-\epsilon_0)|B(r_n)|$ sites belong
to $\D:=\{x: l_n(x)<\delta_0 k_n\}$. Since we suppose 
$|\D|\ge (1-\epsilon_0)|B(r_n)|$, in each {\it good} period, where
$|\RR_i|>2\epsilon_0 |B(r_n)|$,
there are at least $\epsilon_0|B(r_n)|$ sites of $\D$ which are visited. Thus,
$\D$ receives a total of at least $\epsilon_0|B(r_n)|(\delta_0/\epsilon_0)k_n$
 visits.
Necessarily, one site of $\D$ receives more than $\delta_0k_n$ visits, 
and this contradicts the definition of $\D$. Now, from \reff{range-idea} we obtain
\ba{range-logic}
&&P\pare{\acc{\big|\acc{x:\ l_n(x)\ge \delta_0 \frac{n}{|B(r_n)|}}\big|
\ge \epsilon_0|B(r_n)|}\cap \acc{\sigma(r_n)>n}}\cr
&&\qquad\qquad\ge P\pare{\sigma(r_n)>n}-P\pare{
\frac{1}{k_n}\sum_{i=1}^{k_n} X_i> 1-\frac{\delta_0}{\epsilon_0} }.
\ea
Note that by classical estimates $P\pare{\sigma(r_n)>n}\ge 2c_1\exp(-c_2 n/r_n^2)$
for two constants $c_1,c_2$. Finally, the possibility of having $cI(\epsilon_0)$ large,
by reducing $\epsilon_0$, in \reff{range-cheb} allows us to conclude \reff{visits-1}.
\subsection{Proof of the Lower Bound in Region III}\label{sec-III}
We consider $\{X_n>n^\beta\}$. We fix $u=\frac{9}{5}-\frac{6}{5}\beta$, and $v=1-u$.
Note that in Region III, $u$ and $v$ are positive, and
$\zeta_{I\!I\!I\!}=2\beta-2v-u=1-\frac{2}{3}u$. We consider a sequence
of radii with $|B(r_n)|=n^u$ and keep $\epsilon_0$ and $\delta_0$ of Proposition
~\ref{prop-range}. Now, we set $\G:=\{x: l_n(x)>\delta_0 n^v\}$, and use
inequality (2.3) of Lemma 2.1 of \cite{ACa}, since we have assumed that
the $\eta$'s are bell-shaped.
\ba{III-low1}
\PP\pare{\sum_{x\in \Z^d} \eta(x) l_n(x)>n^{\beta} }&\ge &
\PP\pare{ \sum_{\G} \eta(x) \delta_0 n^v >n^\beta }\cr
&\ge & P\pare{ |\G|>\epsilon_0 n^u } P_{\eta}\pare{ \sum_{i=1}^{\epsilon_0 n^u }
\eta_i>\frac{n^{\beta-v}}{\delta_0} },
\ea
where $\{\eta_j,j\in \N\}$ are i.i.d with the same law as $\eta(0)$. Note that
the last probability estimate in \reff{III-low1} on the sum of $\eta$'s is on
the moderate deviations regime, since (i) $\sqrt{ n^u}\ll n^{\beta-v}$,
and (ii) $n^u\gg n^{\beta-v}$. Indeed, (i) is equivalent to $\zeta_{I\!I\!I\!}>0$
which holds, whereas (ii) is equivalent to $\beta< 1$.
Now, in regime (i) and (ii), we have a gaussian lower bound
\be{III-gauss}
P_{\eta}\pare{ \sum_{i=1}^{\epsilon_0 n^u }
\eta_i>\frac{n^{\beta-v}}{\delta_0} }\ge \exp\pare{ -c\frac{n^{2(\beta-v)}}
{\delta_0^2 \epsilon_0 n^u} }=\exp\pare{ -\frac{c}{\delta_0^2 \epsilon_0 }
n^{\zeta_{I\!I\!I\!}}},
\ee
and Proposition~\ref{prop-range} gives the same 
lower bound for $P(|\G|>\epsilon_0 n^u )$.
\section{Upper bounds for deviations estimates for RWRS}
\label{sec-RWRS}
We follow the approach of Section 4 of \cite{ACb}. Thus, we partition
the range of the RW into two domains
$\Dh=\acc{x \in \Z^d: l_n(x) \geq n^{b}}$ and $\Db=
\acc{x \in \Z^d: 0<l_n(x) \leq n^{b}}$, parametrized by a positive $b$. 

According to Section 4 of \cite{ACb}, in each region of interest
we choose $b=\beta-\zeta$, and it is sufficient to find constants $C_1, C_2$
such that for $y$ large enough
\be{low-level}
P\pare{\sum_{x \in \Db} l_n^2(x) > n^{\beta+b} y} \le \exp(-C_1 n^{\zeta})
\ee
and, 
\be{high-level}
P\pare{\sum_{x \in \Dh} l^{\alpha^*}_n(x) \geq n^{\beta-b+\alpha^*b} y}
 \le \exp(-C_2 n^{\zeta})
,\quad\text{where}\quad \alpha^*:=\frac{\alpha}{\alpha-1}.
\ee
\vspace{.5cm}
\noindent{\bf Region I}.
We choose $\beta+b=1$. Since, in Region I, $2\beta-1\le 1/3$,
\reff{low-level} follows from the upper bound in \reff{UB-prop2bis}.
Finally, $b\ge 1/3$ implies that $P(\Dh\not= \emptyset)\le \exp(-Cn^{1/3})$, and
\reff{high-level} holds trivially.

\vspace{.5cm}
\noindent{\bf Region II}.
We choose $b=\beta/(\alpha+1)$. We consider two cases.
\begin{itemize} 
\item First $\beta+b>1$.
The evaluation of $P(\Sigma_n^2>n^{\beta+b}y)$ is straightforward from
the proof of Lemma 2.1 of \cite{ACb} supplied with the moment estimates of the Appendix.
We omit to write this proof, since the argument is by now routine, and
the result reads: for any $\epsilon>0$
\be{old-gamma}
P\pare{ \sum_{x\in \Db} l_n^2(x)> n^{\beta+b} y} \le \exp\pare{
-c n^{\beta+b-\frac{2}{3}-\epsilon}}.
\ee
Now, we can find $\epsilon$ small enough so that in Region II,
$\beta+b-\frac{2}{3}-\epsilon> \beta-b$, which is equivalent to $b>1/3$.
In region II, $b=\beta/(1+\alpha)>1/(4-\alpha)\ge 1/3$. Thus, \reff{low-level} holds.
\item When $\beta+b=1$ (and $\alpha=1$), we have $\zeta_{I\!\!I}=1/3$.
We can take $\epsilon=0$ in \reff{old-gamma},
by \reff{UB-prop2bis}.
\end{itemize}
In order to prove \reff{high-level}, we proceed along the same line
as in \cite{ACb}, and rely on Proposition 3.2 of \cite{ACb}. We omit to repeat
the same computations.

\vspace{.5cm}
\noindent{\bf Region III}.
We choose $5b=\beta+1$. Note that $\beta+b>1$, and with the help of
\reff{old-gamma}, \reff{low-level} follows as soon as $\beta+b-\frac{2}{3}>\beta-b$,
which is equivalent to $\beta>\frac{2}{3}$.

We now prove \reff{high-level}. We consider two cases.
\begin{itemize} 
\item $\alpha\ge d/2$. Condition (0), of Proposition 3.2 of \cite{ACb}, requires
that $\beta-b\le \frac{3}{2}b$ which is equivalent to $\beta\le 1$.
Condition (iii) of the same proposition requires that $\beta<1$.
\item $\alpha< d/2$. We need to check Conditions (i) and (ii) of
Proposition 3.2 of \cite{ACb}. Condition (i) imposes that
\be{cond(i)}
(\beta-b)\pare{ \frac{\alpha^*}{3/2}+1} < \beta-b +\alpha^* b 
\Longleftrightarrow
\beta<\frac{5}{2}b\Longleftrightarrow
\beta<1.
\ee
\reff{cond(i)} is satisfied in Region III.  Condition (ii) requires
\be{cond(ii)}
(\beta-b)\alpha^*<\beta-b+\alpha^*b\Longleftrightarrow
(\beta-b)\frac{\alpha^*-1}{\alpha^*}< b \Longleftrightarrow
\frac{4}{5}\beta -\frac{1}{5}< \alpha(\frac{\beta}{5}+\frac{1}{5})\Longleftrightarrow
\alpha>\frac{4\beta-1}{\beta+1}.
\ee
This last inequality holds in Region III.
\end{itemize}
\section{Appendix}
We have gathered in this section a handy large deviation estimate,
as well as moments computations for variables related to self-intersection times
in dimension 3 and 4.
\subsection{On a large deviation estimate}
\bl{App-UB} Let $\{X,X_1,\dots,X_n\}$ be positive i.i.d. satisfying
\be{eq-UB1}
P(X>u)\le C\exp(-u),\quad\text{with}\quad  C>1.
\ee
We set $\bar X_i=X_i-E[X_i]$, and denote by $c_u=3+e^1+C$. Then, for any 
$\gamma\in ]0,1[$, we have
\be{eq-UB2}
P\pare{\sum_{i=1}^n\bar X_i>x_n}\le \exp\pare{c_un 
\max\pare{\gamma^2 E[X^2],\pare{\gamma^2 E[X^2]}^{1-\gamma}}-
 \frac{\gamma x_n}{2} }.
\ee
\el
\br{app-rem1} Lemma~\ref{App-UB} will serve in regime where $x_n\sim n$.
Estimate \reff{eq-UB2} allows us to take advantage of the smallness of $nE[X^2]/x_n$
to bypass the lack of Cramer's condition. Indeed, assume for instance that instead of
\reff{eq-UB1}, we had for some $\Gamma>0$ (that we think of as a small number
which may depend on $n$) and for $0<\gamma<1$
\be{eq-UB6}
P(X>u)\le C\exp(-\Gamma u),\quad\text{and}\quad 
\max\pare{\Gamma^2 E[X^2],\pare{\Gamma^2 E[X^2]}^{1-\gamma}}\le 
\frac{\Gamma x_n}{4c_un},
\ee
then, the estimate \reff{eq-UB2} would read
\be{eq-UB7}
P\pare{\sum_{i=1}^n\bar X_i>x_n} \le \exp\pare{ -\frac{\gamma\Gamma x_n}{4} }.
\ee
Note that Lemma 1 of \cite{bcr05} does not achieve the same purpose, since
even if $n \Gamma^2 E[X^2]$ were bounded, their proof would yield an
estimate $P(\sum \bar X_i>x_n)\le \exp(-c\Gamma x_n/\log(n))$.
\er
\bpr
Note that for any $\gamma\in ]0,1[$, we use \reff{eq-UB1} and Chebychev to obtain
\be{eq-UB3}
P(X>u)\le \pare{ \frac{E[X^2]}{u^2} }^{1-\gamma} C^{\gamma} e^{-\gamma u}.
\ee
Now, for any $0<\lambda<1$ we decompose $E[\exp(\lambda \bar X)]$ as follows
\ba{eq-UB4}
E[\exp(\lambda \bar X)]&=& E\cro{e^{\lambda \bar X}1\{A\} }
+E\cro{e^{\lambda \bar X}1\{A^c\} }\quad\text{with}\quad
A=\{\lambda X<1\}\cr
&\le & E\cro{e^{\lambda \bar X}\{\lambda \bar X<1\} }+
E\cro{e^{\lambda X}1\{A^c\} }\cr
&\le & E\cro{e^{\lambda \bar X}\{\lambda \bar X<1\} }+e^1P(A^c)+
\int_{1/\lambda}^{\infty}\!\!\!\!\lambda e^{\lambda u} P(X>u) du\cr
&\le & 1+\lambda E\cro{ \bar X1\{\lambda \bar X<1\}}
+2\lambda^2 E[X^2]+e^1P(A^c)+\!\!\int_{1/\lambda}^{\infty}\!\!\!\!
\lambda e^{\lambda u} P(X>u) du\cr
&\le & 1+\lambda E\cro{|\bar X|1\{\lambda \bar X\ge 1\} }+(e^1+2)\lambda^2 E[X^2]
+\lambda\!\!\int_{1/\lambda}^{\infty}\!\!e^{\lambda u} P(X>u) du.
\ea
We have used that for $x\le 1$, $e^x\le 1+x+2x^2$ and that $E[\bar X]=0$. Now,
we choose $2\lambda =\gamma$ and \reff{eq-UB3} to estimate the last term 
in \reff{eq-UB4}
\ba{eq-UB5}
E[\exp(\lambda \bar X)]&=& 1+(3+e^1)\lambda^2E[X^2]+
\lambda C^{\gamma}
\!\!\int_{1/\lambda}^{\infty}\!\!\pare{\frac{E[X^2]}{u^2} }^{1-\gamma} e^{-
\lambda u} du\cr
&\le & 1+(3+e^1)\lambda^2E[X^2]+C^{\gamma}\pare{ \lambda^2 E[X^2]}^{1-\gamma} \cr
&\le & \exp\pare{ c_u \max\pare{\lambda^2 E[X^2],\pare{\lambda^2 E[X^2]}^{1-\gamma}}}.
\ea
The estimate \reff{eq-UB2} follows at once.
\epr
\subsection{Moments computations}
For notational convenience, we keep $n/2$ to denote the integer part
of $n/2$.
\bl{lem-app1} There is $C_0$ such that for $|x|>\sqrt{n}$ and $k< n/2$
\be{app.1}
P_0( S_{n/2-k}=x)\le C_0 P_0(S_{n-k}=x).
\ee
\el
\br{rem-app1} Note that this implies that for $|x|>\sqrt{n}$
\be{app.2}
\sum_{k=0}^n P_0(S_k=x)\le (C_0+1)\sum_{k=n/2}^n P(S_k=x).
\ee
\er
\bpr
Since classical Gaussian estimates gives
\be{app.3}
\frac{C_1 e^{-|x|^2/2k}}{k^{d/2}}\le P_0(S_k=x)\le \frac{C_2 e^{-|x|^2/2k}}{k^{d/2}},
\ee
\reff{app.1} follows if there is a constant $C$, independent of $|x|$ and $k$ such that
\be{app.4}
C \exp\pare{\pare{ \frac{1}{n/2-k}-\frac{1}{n-k} }\frac{|x|^2}{2} }\ge 
\pare{\frac{n-k}{n/2-k} }^{d/2},\quad\text{for}\quad |x|^2>n,\ k<n/2.
\ee
Inequality \reff{app.4} is equivalent to
\be{app.5}
C \exp\pare{\frac{|x|^2}{2(n-k)}\pare{ \frac{n/2}{n/2-k} } }\ge
\pare{ 1+\frac{n/2}{n/2-k} }^{d/2}.
\ee
Thus, since $|x|^2/(n-k)\ge n/(n/2)=1/2$, it is enough to choose 
\[
C:=\sup_{y>1}\pare{ \exp(-\frac{y}{4}) (1+y)^{d/2}}.
\]
We obtain \ref{app.1} by choosing $C_0=C C_2/C_1$.
\epr

We consider $\{\tilde S_n,n\in \N\}$ and independent copy of the random
walk $\{S_n,n\in \N\}$, and denote by $\{\tilde l_n(x), x\in \Z^d\}$ its
local times. Also, we denote $I_n=\sum_{\Z^d} l_n(x) \tilde l_n(x)$.
\bl{lem-app2} In dimension 3, there is a constant $c_3$ such that
$E[I_n]\le c_3 \sqrt{n}$. In dimension 4, there is a constant $c_4$ such that
$E[I_n]\le c_4 \log(n)$.
\el

\br{rem-app2} Note that when $n=2^N$, and $\{ I_k^{(l)},\ k=1,\dots,2^l\}$ are
independent copies with the same distribution as $I_{2^{N-l}}$, we have both for $d=3,4$
constants $C_3$ and $C_4$ such that
\be{app.6}
E\cro{ I_k^{(l)} }\le
\left\{ \begin{array}{ll}
                C_3  \sqrt{2^{N-l}} & \mbox{ for } d=3 \, 
                \\
                 C_4  (N-l) & \mbox{ for } d \geq 4 \, ,
                \end{array}
        \right.
\quad\text{and}\quad
E\cro{ \sum_{l=1}^{N-1}\sum_{k=1}^{2^{l-1}} I_k^{(l)} }\le \frac{C_d}{\sqrt{2}-1} 2^N.
\ee
\er
\bpr
If we denote by $\gamma_d$ the probability of not returning to 0, i.e. 
$P_0(H_0=\infty)=\gamma_d>0$, then $E_0[l_{\infty}(0)]=1/\gamma_d$, and
\be{app.7}
E[I_n]=\sum_{x\in \Z^d} \pare{ E_0[l_n(x)]}^2\le
\frac{1}{\gamma_d^2}\sum P_0(H_x\le n)^2\le R_{n,1}+R_{n,2},
\ee
with 
\ba{app.8}
R_{n,1}:=\frac{1}{\gamma_d^2}\sum_{|x|\le \sqrt{n}} P_0(H_x\le n)^2&\le&
C \sum_{|x|\le \sqrt{n}} \frac{1}{1+|x|^{2(d-2)}}\cr
&\le&C' \int_1^{\sqrt{n}} \frac{x^{d-1}}{x^{2(d-2)}} dx\le C'' \left\{ \begin{array}{ll}
                \sqrt{n} & \mbox{ for } d=3 \, ,
                \\
                \log(\sqrt{n} )  & \mbox{ for } d \geq 4 \, .
                \end{array}
        \right.
\ea
Now, for $R_{n,2}$, we note that $P_0(H_x\le n)\le P_0(S_0=x)+\dots+P_0(S_n=x)$,
and use Lemma~\ref{lem-app1}
\be{app.9}
R_{n,2}:=\frac{1}{\gamma_d^2}\sum_{|x|>\sqrt{n}} P_0(H_x\le n)^2\le
\frac{(C_0+1)^2}{\gamma_d^2}\sum_{|x|>\sqrt{n}} \pare{ \sum_{k=n/2}^n P_0(S_k=x) }^2.
\ee
Now, note that from \reff{app.3}, there is $C$ such that for $n\ge k\ge n/2$
\be{app.10}
P_0(S_k=x)\le \frac{C e^{-x^2/(2n)}}{n^{d/2}}.
\ee
Thus,
\be{app.11}
\sum_{|x|>\sqrt{n}} \pare{ \sum_{k=n/2}^n P_0(S_k=x) }^2\le 
\sum_{|x|>\sqrt{n}} \frac{C^2 e^{-x^2/n}}{n^{d-2}}.
\ee
Finally, there is a constant $C'$ such that
\be{app.12}
R_{n,2}\le C'\int_{\sqrt{n}}^n\!\! \frac{e^{-x^2/n}}{n^{d-2}} x^{d-1} dx\le
C' n^{2-d/2} \int_1^{\infty} e^{-u^2} u^{d-1} du.
\ee
The result follows as we gather \reff{app.9} and \reff{app.12}.
\epr

We denote now $\D_n(z):=\{x:\ l_n(x)>z\}$. The following Lemma 
estimates the first two moments of $\tilde l_n\pare{\D_n(z) }$.
\bl{lem-app3}
There are positive constants $\kappa_3,\kappa_4,C_3,C_4$ such that
\be{app.30}
E\otimes \tilde E\cro{ \tilde l_n\pare{\D_n(z) } }\le C_d\left\{ \begin{array}{ll}
                n^{2/3}\exp(-\frac{2}{3}\kappa_3 z) & \mbox{ for } d=3 \, ,
                \\
                \sqrt{n} \exp(-\frac{\kappa_4}{2} z)  & \mbox{ for } d \geq 4 \, .
                \end{array}
        \right.
\ee
Moreover, we also have constants $C_3'$ and $C_4'$ such that
\be{app.13}
E\otimes \tilde E\cro{ \tilde l_n\pare{\D_n(z) }^2 }\le C_d'\left\{ \begin{array}{ll}
                n^{4/3}\exp(-\kappa_3 z) & \mbox{ for } d=3 \, ,
                \\
                n \exp(-\kappa_4 z)  & \mbox{ for } d \geq 4 \, .
                \end{array}
        \right.
\ee
\el
\bpr
We have seen in \cite{ACa} that when $d\ge 3$,
there is $c_d$ independent of $n$ and of the domain $\Lambda$ such that
\be{app.29}
\sup_x E_x\cro{ l_n\pare{\Lambda} }\le c_d |\Lambda|^{2/d}.
\ee
Thus, using Holder's inequality
\be{app.31}
E\otimes\tilde E\cro{\tilde l_n\pare{\D_n(z)} }\le  c_d E\cro{ |\D_n(z)|^{2/d}}
\le c_d \pare{ E\cro{ |\D_n(z)|} }^{2/d}.
\ee
Now, since the expected number of visited sites at time $n$, is of order $n$,
we have
\be{app.18}
|\D_n(z)|=\sum_{x\in \Z^d} 1\{l_n(x)>z\},\quad\text{and}\quad
\sup_x E_x\cro{|\D_n(z)|}\le c_d' n e^{-\kappa_d z}.
\ee
Thus,
\be{app.32}
E\otimes\tilde E\cro{\tilde l_n\pare{\D_n(z)} }\le  c_d\pare{c_d' n 
e^{-\kappa_d z}}^{2/d}.
\ee
Inequality \reff{app.30} follows at once.
We now prove \reff{app.13}. First note that
\be{app.14}
 \tilde l_n\pare{\D_n(z)}^2=2\sum_{x,y\in \D_n(z),}\sum_{k<k'\le n} 1\{
\tilde S_k=x, \tilde S_{k'}=y \}+\quad \tilde l_n\pare{\D_n(z)}.
\ee
Now, we average only over the walk $\{\tilde S_n\}$
\ba{app.15}
\tilde E\cro{\tilde l_n\pare{\D_n(z)}^2}&=& 2\sum_{x,y\in \D_n(z)}\sum_{k<k'\le n}
\tilde \P_0(\tilde S_k=x)\tilde P_x(\tilde S_{k'-k}=y)+
\tilde E\cro{ \tilde l_n\pare{\D_n(z)} }\cr
&\le& 2\sum_{x\in \D_n(z)} \sum_{k\le n} \tilde P_0(\tilde S_k=x)\tilde
E_x\cro{ \tilde l_n\pare{\D_n(z)} } +\tilde E\cro{ \tilde l_n\pare{\D_n(z)} }\cr
&\le & 2\pare{\sup_x \tilde E_x\cro{ \tilde l_n\pare{\D_n(z)} } }^2
+\tilde E\cro{ \tilde l_n\pare{\D_n(z)} }.
\ea
From \reff{app.29} we obtain
\be{app.16}
\pare{\sup_x \tilde E\cro{\tilde l_n\pare{\D_n(z)}}}^2\le C_d^2 |\D_n(z)|^{4/d}.
\ee
We average now with respect to the random walk $\{S_n\}$, (and use Jensen's inequality
in $d=3$)
\be{app.17}
E\cro{\pare{\sup_x \tilde E\cro{\tilde l_n\pare{\D_n(z)}}}^2}
\le C_d^2\left\{ \begin{array}{ll}
                \pare{E\cro{|\D_n(z)|^2}}^{2/3} & \mbox{ for } d=3 \, ,
                \\
                E\cro{ |\D_n(z)| } & \mbox{ for } d \geq 4 \, .
                \end{array}
        \right.
\ee
Finally, note that
\be{app.19}
|\D_n(z)|^2\le |\D_n(z)|+2\sum_{x\not= y} 1\{H_x<H_y\le n,\ l_n(y)>z\}.
\ee
Taking the expectation in \reff{app.19}, we obtain
\ba{app.20}
E\cro{ |\D_n(z)|^2 } 
&\le & E\cro{ |\D_n(z)|}+ 2\sum_{x}E\cro{ 1\{H_x< n\}\sum_y P_x\pare{ l_n(y)>z}}\cr
&\le & 
\sup_x E_x\cro{ |\D_n(z)|}\pare{1+2 E_0\cro{\{x: l_n(x)>0\} }}\le C n^2 e^{-\kappa_d z}.
\ea
This concludes the proof.
\epr


\begin{thebibliography}{99}
\bibitem{ACa} Asselah, A., Castell F.,
{\em A note on random walk in random scenery.  } 
To appear in Annales de l'I.H.P., also arXiv:math.PR/0501068.

\bibitem{ACb} Asselah, A., Castell F., 
{\em Self-Intersection Times for Random Walk, and Random Walk in Random 
Scenery in dimensions $d\ge 5$.  } Preprint 2005, arXiv:math.PR/0509721 

\bibitem{bc04} Bass, R.F., Chen, X., 
{\em Self-intersection local time: critical exponent, large deviations, 
and laws of the iterated logarithm.}
Ann. Probab.  32  (2004),  no. 4, 3221--3247.

\bibitem{bcr05} Bass, R.F., Chen, X., Rosen, J., 
{\em Large deviations for renormalized self-intersection local times of stable 
processes.}  Ann. Probab.  33  (2005),  no. 3, 984--1013.

\bibitem{bcr05b} Bass R.F., Chen X., Rosen J.
{\em    Moderate deviations and laws of the iterated logarithm 
for the renormalized self-intersection local times of planar random walks}
Preprint 2005, arXiv, math.PR/0506414. 

\bibitem{bcr06} Bass R.F., Chen X., Rosen J.
{\em  Moderate deviations for the range of planar random walks}
Preprint 2006, arXiv, math.PR/0602001 .

\bibitem{bolt}
van den Berg, M.; Bolthausen, E.; den Hollander, F. 
{\em Moderate deviations for the volume of the Wiener sausage.}
Ann. of Math. (2)  153  (2001),  no. 2, 355--406. 

\bibitem{chen-li}
Chen, Xia; Li, Wenbo V.
{\em Large and moderate deviations for intersection local times.}
Probab. Theory Related Fields  128  (2004),  no. 2, 213--254.

\bibitem{DV} Donsker, M. D.; Varadhan, S. R. S. 
{\it Asymptotic evaluation of certain Markov process for large
time. I. II. III. IV.}
Comm. Pure Appl. Math. {\bf 28} (1975), 1--47; ibid {\bf 28} (1975), 279--301; 
ibid {\bf 29} (1976), no 4, 389--461; ibid {\bf 36} (1983), no 2, 183--212.

\bibitem{HGK} Gantert, N.; van der Hofstad, R.; K\"onig, W.
{\em Deviations of a random walk in a random scenery with
stretched exponential tails}. Preprint 2004.
arXiv:math.PR/0411361.

\bibitem{GKS} Gantert, N.; K\"onig, W.; Shi, Z.
{\em Annealed deviations of random walk in random scenery}
Preprint 2004.  arXiv.:math.PR/0408327.


\bibitem{legall} Le Gall, J.F.;
{\em Sur le temps local d'intersection du mouvement brownien plan et la m\'ethode
de renormalisation de Varadhan.}
S\'eminaire de probabilit\'es, XIX, 1983/84,  314--331,
Lecture Notes in Math., 1123, Springer, Berlin, 1985.

\bibitem{Ma91} Mansmann, U.;
{\it The free energy of the Dirac polaron, an explicit solution.}
Stochastics Stochastics Rep.  34  (1991),  no. 1-2, 93--125.


\end{thebibliography}
\end{document}